\nonstopmode
\documentclass[12pt]{article}

\usepackage{anysize}
\usepackage{amsmath,amsthm, amssymb}
\usepackage{amssymb}
\usepackage{latexsym}
\usepackage{url,paralist}
\usepackage{color}

\usepackage{bbm}

\newcommand{\sr}{\overline {\overline r}}
\newcommand{\halffl}[1]{\lfloor #1/2 \rfloor}
\renewcommand{\vec}[1]{\boldsymbol #1}
\newcommand{\Prob}[1]{\mbox{\normalfont Prob}\left[ #1 \right]}
\newcommand{\so}{o}
\newcommand{\bo}{O}

\newcommand{\mat}{{\cal M}_{d}}
\newcommand{\sing}{{\cal S}}
\newcommand{\reg}{{\cal R}}
\newcommand{\gd}{{\cal G}}

\newcommand{\pd}{P_s(d)}
\newcommand{\ed}{E(d)}

\newcommand{\naturals}{\mathbbm{N}}

\newcommand{\reals}{\mathbbm{R}}

\newcommand{\suchthat}{\ | \ }

\theoremstyle{plain}
\newtheorem{lemma}{Lemma}[section]
\newtheorem{theorem}[lemma]{Theorem}
\newtheorem{korollar}[lemma]{Corollary}
\theoremstyle{definition}

\theoremstyle{remark}

\begin{document}

\title{Singular $0/1$-matrices, and the\\
hyperplanes spanned by random $0/1$-vectors}
\author{Thomas Voigt\thanks{%
Supported by Deutsche Forschungs-Gemeinschaft (DFG) ZI 475/3}
\setcounter{footnote}{6}
\qquad G\"unter M. Ziegler\thanks{%
Partially supported by Deutsche Forschungs-Gemeinschaft (DFG), FZT86 and ZI 475/4}\\
\small Inst.\ Mathematics, MA 6-2\\[-1.6mm]
\small TU Berlin, D-10623 Berlin, Germany\\[-1.6mm]
\small\url{{tvoigt,ziegler}@math.tu-berlin.de}}
\date{}

\maketitle

\begin{abstract}\noindent
Let $\pd$ be the probability that a random 
$0/1$-matrix of size $d \times d$ is 
singular, and let $\ed$ be the expected number of $0/1$-vectors in 
the linear subspace spanned by $d-1$ random independent $0/1$-vectors.
(So $\ed$ is the expected number of cube vertices 
on a random affine hyperplane spanned by vertices of the cube.)

We prove that bounds on $\pd$ are equivalent to
bounds on $\ed$: 
$\pd = \left(2^{-d} \ed + \frac{d^2}{2^{d+1}} \right) (1 + \so(1))$.

We also report about computational experiments pertaining to these numbers.
\end{abstract}

\section{Introduction}

$0/1$-polytopes arise naturally in a great variety of interesting
contexts, including a prominent role in combinatorial optimization, yet
some basic characteristics of ``typical'' (that is, random)
$0/1$-polytopes are unknown. (For a survey of a variety of 
aspects of $0/1$-polytopes see \cite{Zie99}.) 

One of the key open questions in this context is rather notorious: 
\begin{itemize}
\item
Pick $d+1$ random
vertices of the $d$-cube independently (with respect to the uniform
distribution). 
What is the probability that these vectors do not form a $d$-simplex?
\end{itemize}
If we assume without loss of generality that one of these points
is the origin $\vec{0}$ the question can be rephrased: Let
$C^d=[0,1]^d$ be the $d$-dimensional unit hypercube, and let
\[ \mat\ \ :=\ \ \{0,1\}^{d \times d} \]
be the set of all $0/1$-matrices of size~$d\times d$.
\begin{itemize}
\item
What is the asymptotic behaviour of the probability
\[ 
\pd\ \ :=\ \ \Prob{\det(M) = 0 \suchthat M \in \mat} \]
that a random square $d$-dimensional $0/1$-matrix
is singular?
\end{itemize}
This central but difficult question has received careful
attention; see Koml\'os \cite{Kom67}, Bollob\'as \cite{Bol85}, 
Kahn, Koml\'os \& Szemer\'edi \cite{KKS95}.
It has been conjectured that 
\[
\pd\ \ =\ \ \frac{d^2}{2^d}(1 + \so(1)), 
\]
which is essentially the probability
that two rows or two columns of a random matrix are equal. However,
the known upper bounds are far off this mark; currently the best upper
bound is $\pd < (1-\varepsilon)^d$, for some rather small
$\varepsilon>0$.
(This was proved by Kahn, Koml\'os and Szemer\'edi in \cite{KKS95} with $\varepsilon=0.001$.)
\smallskip

A closely related problem is as follows:
\begin{itemize}
\item
Given $r$ random vertices $\vec{v}_1, \ldots, \vec{v}_r$ of $C^d$, what is 
the expected number of $0/1$-vectors in the affine subspace spanned by
these vectors? 
\end{itemize}
Improving a result by Odlyzko \cite{Odl88}, Kahn, Koml\'os \& 
Szemer\'edi derived in \cite{KKS95} that there exists a
constant $C$ independent from $d$ such that the probability that such
an affine subspace contains any $0/1$-vector
other than $\vec{v}_1, \ldots, \vec{v}_r$ is 
$4 \binom{r}{3}(\frac{3}{4})^d (1 + \so(1))$,
provided that $r < d-C$. 
However, so far no results were known for the case $r=d$. 

In this paper we will show that determining the expected number of
vertices of $C^d$ in the affine
subspace spanned by $d$ random vertices of $C^d$ is just
as hard as determining $\pd$. More precisely, let
$\gd$ denote the set of all linearly independent $(d-1)$-sets of $0/1$-vectors of
length $d$ and for a set $S$ of arbitrary vectors let $v(S)$ be the number of $0/1$-vectors 
in the linear subspace spanned by $S$.
Then the following theorem holds.

\begin{theorem}
\label{th1}
Let 
\[ 
\ed\ \ :=\ \ \frac{1}{|\gd|} \sum_{G \in \gd} v(G) 
\]
be the expected number of $0/1$-points on the hyperplane spanned 
by a random linearly independent set of $d-1$ $0/1$-vectors.
Then
\[ \pd = \left(\frac{1}{2^d} \ed + \frac{d^2}{2^{d+1}} \right) (1 + \so(1)).
\] 
\end{theorem}

We can give a (trivial) lower bound for $\ed$ by just considering the
${d \choose 2} + d$ ``fat'' hyperplanes (faces $x_i = 0$ and
hyperplanes $x_i - x_j = 0$) containing $2^{d-1}$ vertices each. 
Since $d-1$  points chosen randomly from such a hyperplane 
span the hyperplane with probability $1-(1 - \varepsilon)^{d-1}$ (according to \cite{KKS95}) 
it is easy to verify that $\ed \geq \frac{d^2}{2}(1 + \so(1))$.

In fact the conjectured upper bounds on $\pd$ and $\ed$
are strictly equivalent:
\begin{korollar}
\label{kor2}
As $d \rightarrow \infty$, 
\begin{eqnarray*}
\pd &=& \frac{d^2}{2^d}(1 + \so(1))\\
\noalign{if and only if}
\ed &=& \frac{d^2}{2}(1 + \so(1)).
\end{eqnarray*}
\end{korollar}

Using symmetry we could switch to an affine version, replacing $\gd$ by the
set of affinely independent $d$-sets of $0/1$-vectors and checking the
expected value of $0/1$-vectors in a hyperplane spanned by such a set.
However, for the purpose of this paper the linear version will be
more convenient to handle; so we will consider only hyperplanes
containing the origin~$\vec{0}$.

To our knowledge the problem of determining the expected number of
$0/1$-vectors on a hyperplane $h$ spanned by random vertices of $C^d$ 
has not been studied independently yet. Some basic results were derived in 
\cite{Bol85} and \cite{KKS95} by  
examining the structure of the defining equations $\vec{a}$ for planes
$h = \{\vec{x} \in \reals^d \suchthat \vec{a}^t \vec{x} = 0\}$ (which
is perhaps the most natural approach). 
The lemma of Littlewood-Offord (see Section~\ref{sec:preparation})
a classical tool:  
It states that if all $a_j$ are 
nonzero then the number of $0/1$-points in this plane is at 
most~$\binom{d}{\halffl{d}}$. If the coefficients satisfy additional
conditions, this number can be reduced considerably  
(see Hal\'asz \cite{Hal75} \cite{Hal77}). In order to obtain such conditions
it would be of considerable interest to learn more about the distribution
of determinants of $0/1$-matrices: If $d-1$ vectors span a
hyperplane and we write these vectors into a $d \times(d-1)$ matrix $M$, 
then a defining equation $\vec{a}^t \vec{x}=0$ is given by 
$a_j = (-1)^j \det(r_j(M))$, where $r_j(M)$ is the matrix obtained from $M$
by deleting the $j$-th row. \\

The rest of this paper is organized as follows: In
Section~\ref{sec:preparation} we state some consequences 
of the Littlewood-Offord lemma. The proof of Theorem \ref{th1} 
is given in Section~\ref{sec:proof}. In Section~\ref{sec:experiments}
we present some experimental estimates of $\pd$ for~$d \leq 30$. 
\medskip

{\bf Some definitions.} \\
We use standard vector notation $\vec{a} = (a_1, \ldots, a_d)^t$, where $d$ denotes 
the dimension. The expected value of a random variable $X$ 
is denoted by $E[X]$; the probability of an event $Y$ is $\Prob{Y}$. 
Define $r(F)$ as the (linear) rank of a family or set of vectors $F$.

The next definition is useful for partitioning sets of matrices into
subsets with ``nice'' properties and was frequently used in the
analysis of $0/1$-matrices (see \cite{Bol85} or \cite{KKS95}).
Given a $d \times d$ matrix $M$ we 
define the \emph{strong rank} $\sr(M)$ as 
the largest $k\le d$ such that all
$k$-subsets of columns from $M$ are independent.
(Equivalently, it is the largest $k$ such that the
truncation to rank~$k$ of the matroid given by the columns of the matrix
$m$ is uniform of rank~$k$.) We also consider
the strong rank of sets and of families of $d$-dimensional vectors.

\section{The Littlewood-Offord lemma}
\label{sec:preparation}

The ``Littlewood-Offord lemma'' is a
classical tool  \cite{Bol85} \cite{Odl88} for obtaining upper bounds on~$\pd$.

\begin{lemma}[Littlewood-Offord]
\label{LO2}
Let $s \in \reals$, $n \in \naturals$ and let $a_i \in \reals$ with
$|a_i| \ge1$ for $1\le i\le n$.
Then at most $\binom{n}{ \lfloor \frac{n}{2} \rfloor}$
of the $2^n$ sums $\sum_{i=1}^{n} \varepsilon_i a_i$, 
$\varepsilon_i = \pm 1$ fall in the open interval $(s-1, s+1)$.
\end{lemma}

\begin{korollar}
\label{LO}
Let $a_i \in \reals, i=1, \ldots ,n$ with at least
$t$ of the  $a_i$ nonzero. Then at most
$\binom{t}{\lfloor t/2\rfloor}2^{n-t}\approx\frac{2^n}{\sqrt{\frac{\pi}{2}t}}$ 
of the $2^n$ sums
$\sum_{i=1}^{n} \varepsilon_i a_i$, $\varepsilon_i \in \{0, 1\}$ can have
the same value.
\end{korollar}

As observed in \cite{KKS95}, this lemma suffices to show that 
with very high probability the strong rank of a 
random $0/1$-matrix is either close to~$d$ or at most~$1$.

\begin{lemma}
\label{LOLemma}
Let $M \in \mat$ be a random matrix. Let $E$ be the event that 
$M$ has a $d \times (k+1)$ submatrix of strong rank $k$
for some $k \in \{2, \ldots, d - 3 \frac{d}{\ln(d)} \}$.
Then for large~$d$,
\[ 
\Prob{ E }\ \ \le\ \ 2^{-d}.
\]
\end{lemma}

\begin{proof}
The proof follows \cite[Chapter 14.2]{Bol85} (see also \cite[Section 3.1]{KKS95})
and is sketched here for the reader's convenience.

Let $M$ be a random $0/1$-matrix and $k < d$. 
If $M$ contains
$k+1$ columns $\vec{c}_1, \ldots, \vec{c}_{k+1}$ of strong rank
$k$ then clearly we can find a $k \times (k+1)$ submatrix
of $M$ of strong rank $k$ by deleting $d-k$ linearly dependent rows
from $(\vec{c}_1, \ldots, \vec{c}_{k+1})$. 

If we want to upper bound the probability that arbitrarily
chosen columns $\vec{c}_1, \ldots, \vec{c}_{k+1}$ have strong
rank $k$, then it suffices to give an upper bound on the probability 
that $\vec{c}_1, \ldots, \vec{c}_{k+1}$ have rank $k$
conditioned on the event that an arbitrary $k \times (k+1)$ submatrix 
$\tilde M$ of $(\vec{c}_1, \ldots, \vec{c}_{k+1})$ has strong rank $k$:

$\tilde M$ has strong rank $k$ if and only if the last column
of $\tilde M$ is a unique linear combination of the first $k$ columns
and all coefficients in this combination are non-zero. 
Under this condition the probability that any of the
remaining $d-k$ rows of $\vec{c}_1, \ldots, \vec{c}_{k+1}$ satisfy the linear dependency
equation defined by $\tilde M$ is at most 
$2^{-k} \binom{k}{\lfloor\frac{k}{2} \rfloor}$ by Lemma~\ref{LO},
so the probability that $\vec{c}_1, \ldots, \vec{c}_{k+1}$
have rank $k$ is at most $(2^{-k}
\binom{k}{\lfloor\frac{k}{2} \rfloor})^{d-k}$.
Since there are at most ${d \choose k} {d \choose k+1}$ such submatrices
$\tilde M$ we find

\[ 
\Prob{\sr(M) = k}\ \ \leq\ \  \binom{d}{k} \binom{d}{k+1} 
\left(2^{-k} \binom{k}{\lfloor \frac{k}{2} \rfloor} \right) ^{d-k}. \]

We derive
\begin{equation}
\label{eq:losum}
 \sum_{k=3}^{\lfloor d - 3 \frac{d}{\ln(d)} \rfloor}  \binom{d}{k} \binom{d}{k+1} \left(2^{-k}
\binom{k}{\lfloor \frac{k}{2} \rfloor} \right) ^{d-k}\ \ \le\ \ 2^{-d}
\end{equation}
by checking that each summand in~(\ref{eq:losum}) is at most
$\frac{1}{d 2^d}$ if $d$ is large (using Stirling's formula and elementary, but
somewhat tedious calculations). 

To complete the proof of Lemma~\ref{LOLemma} we observe that the event
$\sr(M) = 2$ depends on the existence of three columns $\vec{m}_i, \vec{m}_j, \vec{m}_k$ such
that $\vec{m}_i + \vec{m}_j = \vec{m}_k$,
which happens with probability $\Theta(d^3 (\frac{3}{8})^d)$.
\end{proof}

\begin{korollar}
\label{LOKorollar}
Let $M \in \mat$ be a random matrix.
Then 
\[ 
\Prob{\sr(M) \le d - 3 \frac{d}{\ln(d)}}
\ \ \le\ \ 
 \frac{d^2}{2^{d+1}}(1 + \mbox{\normalfont o}(1)).
\]
\end{korollar}

\newpage
\section{Proof of Theorem \ref{th1}}
\label{sec:proof}

Let $\sing \subset \mat$ be the set of singular matrices and $\reg = \mat \setminus \sing$.
We will partition $\sing$ into subsets 
$\sing_j \subset \sing, j \in \{1, \ldots, 4\}$ and
derive precise bounds on the sizes of two of these sets 
in terms of $|\gd|$ and $\ed$.
The other two sets are small. This allows us to
estimate the value $\pd = \frac{|\sing|}{|\sing|+|\reg|}$.

Let  $N_d :=\lfloor d - \frac{3d}{\ln(d)}\rfloor$ and partition $\sing$ into the disjoint sets
\begin{eqnarray*}
\sing_1 & := & \{ M \in \mat  \suchthat r(M) = d-1, \ \sr(M) = 1 \} \\
\sing_2 & := & \{ M \in \mat  \suchthat r(M) = d-1, \ \sr(M) > N_d \} \\
\sing_3 & := & \{ M \in \mat \suchthat \sr(M) \in \{0, 2, \ldots, N_d\} \} \\
\sing_4 & := & \{ M \in \mat  \suchthat r(M) < d-1, \ \sr(M) = 1 \mbox{ or } 
\sr(M) > N_d \}. 
\end{eqnarray*}
We will give precise estimates for the sizes of 
the sets $\reg$, $\sing_1$, and $\sing_2$,
and check that the sets $\sing_3$ and $\sing_4$ are small enough. 
More precisely, we will show that
\begin{eqnarray}
\label{reg_est}
| \reg   | & = &  |\gd|\, d!  \frac{2^d - \ed}{d}  \\
\label{sing1_est}
| \sing_1 | & = &  |\gd|\, d!  \frac{d-1}{2} \\
\label{sing2_est}
| \sing_2 | & = &  |\gd|\, d! \frac{\ed}{d}  (1 + \so(1)) \\
\label{sing1_sing2_est}
|\sing_1| & \leq & |\sing_2|(1 + \so(1)) \\
\label{sing3_est}
| \sing_3 | & \leq &  \frac{c_1}{d} |\sing_1| \\
\label{sing4_est}
| \sing_4 | & \leq &  \frac{c_2}{\sqrt{d}}  (|\sing_1| + |\sing_2|) 
\end{eqnarray}
for some constants $c_1, c_2>0$. \\

\begin{itemize}
\item While most matrices from $\mat$ with two equal columns are in $\sing_1$, most matrices with 
two equal rows lie in $\sing_2$. To see this,
pick a random $(d-1) \times d$ matrix $N = (\vec{n}_1, \ldots, \vec{n}_d)$. 
Using the result of Kahn, Koml\'os and Szemer\'edi \cite{KKS95}
 that $\pd \leq (1-\varepsilon)^d$ for some $\varepsilon \geq 0.001$,
we obtain $d (1-\varepsilon)^{d-1}$ as an upper bound on the 
probability that at least one of the
$(d-1) \times (d-1)$ submatrices $c_j(N)$ is singular,
where $c_j(N)$ is the matrix obtained
from $N$ by deleting the $j$-th column $\vec{n}_j$.  
Cramer's rule gives $\sum_{j=1}^d (-1)^j d_j \vec{n}_j = \vec{0}$ for
the determinants $d_j = \det(c_j(N))$.
Thus, $N$ has strong rank $d-1$ if all determinants are nonzero, 
which establishes (\ref{sing1_sing2_est}):
\[ |\sing_1|  \leq  |\sing_2|(1 + \so(1)) \]

\item By Lemma~\ref{LOLemma} a random matrix $M \in \mat$ lies in $\sing_3$ with
probability at most $(d+1) 2^{-d}$. The probability that two columns are equal 
is $d^2 2^{-d-1} (1+\so(1))$. Again almost all matrices with two identical 
columns have strong rank $d-1$ and are in $\sing_1$ (up to an exponentially small subset),
which implies (\ref{sing3_est}):
\[|\sing_3| = \bo(\frac{1}{d} |\sing_1|) \]
\end{itemize}

For each matrix $M \in \reg \cup \sing_1 \cup \sing_2$
there is at least one $G \in \gd$ that is a subset of the column set of $M$. 
The estimates (\ref{reg_est}), (\ref{sing1_est}) and (\ref{sing2_est}) are obtained by examining this in detail:

\begin{itemize}
\item For each $G \in \gd$ we have exactly $\frac{d!}{2} (d-1)$ matrices 
from $\sing_1$ containing only columns from $G$ (since we have $d-1$ choices for a duplicate
column and $\frac{d!}{2}$ permutations). This gives (\ref{sing1_est}):
\[ | \sing_1 | \ \ =\ \ |\gd| \frac{d!}{2} (d-1) \]

\item For any fixed $G \in \gd$ we can construct $d! (v(G)-d)$ 
different matrices $S \in \sing_2 \cup \sing_3$
(using columns from $G$ and an additional nonzero column in the span
of $G$ that is not in $G$).
Summing over $G \in \gd$ we obtain 
$d! \ed (1 + \so(1))|\gd|$
matrices in $\sing_2$, since (\ref{sing1_sing2_est}) and (\ref{sing3_est}) imply
that $|\sing_3|$ is small compared to $|\sing_2|$.
 On the other hand each matrix $M \in \sing_2$ is
constructed $\sr(M)+1$ times: If $M \in S_2$ and $M \vec{a} = 0$ 
for some $\vec{a} \neq \vec{0}$ then $|\mbox{\normalfont supp}(\vec{a})| = \sr(\color{red}{M})+1$ (equality holds
since $r(M)=d-1$)  and
$\{\vec{m}_1, \ldots,\vec{m}_{k-1},\vec{m}_{k+1}, \ldots,\vec{m}_d\}$ 
is independent if and only if
$a_k \neq 0$. 
This gives (\ref{sing2_est}):
\begin{eqnarray*}
|\sing_2| & = & \frac{1}{d - \so(d)}d! \ed (1+\so(1))|\gd| \\
&=& d! \frac{\ed}{d}  (1+\so(1)) |\gd|. 
\end{eqnarray*}

\item Similarly, we get $d!(2^d - \ed)|\gd|$ 
matrices in $\reg$ and each matrix $M \in \reg$ is constructed  $d$ times.
This gives (\ref{reg_est}):
\[ |\reg|  =  \frac{d!}{d}\left( 2^d - \ed \right) |\gd|. \]
\end{itemize}

A little more work is required for the upper bound (\ref{sing4_est}) on $|\sing_4|$. So far we established
an upper bound on the number of matrices of rank $d-1$ in terms of the number
of regular matrices. A similar argument will be used to show that 
for any $k \leq d-2$ there are significantly 
fewer matrices of rank $k$ than matrices of 
rank $k+1$, which gives the desired result:

\begin{enumerate}
\item 
First consider the matrices $\hat \sing$ with the property
that the rows or the columns admit more than one
trivial dependency (i.e. zero-vectors or pairs of identical vectors). This probability
is dominated by the probability that a matrix has two pairs of identical
rows or columns, which happens with probability $\bo({d \choose 4} 2^{-2d})$,
so clearly $|\sing_4 \cap \hat\sing|$ is exponentially smaller than $\frac{1}{\sqrt{d}}(|\sing_1| + |\sing_2|)$. 

\item  
Let $\check \sing$ be the set of matrices whose columns or rows
have a subset with strong rank in $\{2, \ldots, N_d\}$.
Lemma~\ref{LOLemma} gives that this happens with probability
of at most $2^{-d}$, while the probability that two columns are equal 
is $d^2 2^{-d-1} (1+\so(1))$. This implies
$|\check \sing| \leq \bo(\frac{1}{d^2} |\sing_1|)$. 

\item To estimate the number of the remaining matrices in $\sing_4$, we use similar
techniques as in \cite[Chapter 14.2]{Bol85}: \\

We can use the Littlewood-Offord lemma to give an upper bound on the number of
$0/1$-vectors in the span of a set of vectors ${\cal C}$:
Let $\vec{a}$ be in the orthogonal space
of ${\cal C}$, i.e. $\vec{a}^t \vec{c} = 0$ for all $\vec{c} \in {\cal C}$. 
Clearly  all vectors $\vec{v}$ in the span of ${\cal C}$ satisfy $\vec{a}^t \vec{v} = 0$.
If $s$ is the number of nonzero entries in $\vec{a}$
then Lemma \ref{LO}
assures us that the span of ${\cal C}$ contains at most 
${s \choose \halffl{s}} 2^{d-s}$ $0/1$-vectors. \\

Let $\sing_4(k)$ be the matrices in $\sing_4 \setminus (\hat\sing \cup \check\sing)$  
of rank $k$. For a fixed $k \leq d-2$ and $m \in \sing_4(k)$ 
we know that the columns and rows of $m$ 
 admit at most one trivial dependency (by excluding $\hat\sing$)
and that neither rows nor columns have a submatrix of strong rank 
between $2$ and $N_d$ (by excluding $\check\sing$). 
Thus both $\ker(m)$ and $\ker(m^t)$ 
contain vectors with more than $N_d$ nonzero
entries, since they are are at least $2$-dimensional. 
Choose any such vectors $\vec{a} \in \ker(m)$ and $\vec{b} \in \ker(m^t)$.

If $m$ is chosen uniformly at random from $\sing_4(k)$, then
the probability that $a_d \neq 0$  is at least $\frac{N_d}{d} = 1 - \frac{3}{\log d}$.
If we condition on this event (that the last column of $m$ is a nontrivial
linear combination of the remaining columns) and consider all  $0/1$-matrices having
the same first $d-1$ columns as $m$, then (by the observation above)
at most $ 2^{d-N_d} {N_d \choose \halffl{N_d}}$ of these matrices 
have rank $k$, since the last column $\vec{v}$ has to satisfy 
$\vec{b}^t \vec{v} = 0$. 
Stirling's formula implies that 
$ 2^{d-N_d} {N_d \choose \halffl{N_d}} \approx 2^d \sqrt{\frac{2}{\pi N_d}} = \bo(\frac{1}{\sqrt{d}} 2^d)$.

Removing the condition $a_d \neq 0$ changes 
the number of matrices only by a factor of $1 + \frac{3}{\log d}$,
so we find that 
\[ 
|\sing_4(k)| = \left\{ 
\begin{array}{lcl}
\bo( \frac{1}{\sqrt{N_d}}|\sing_4(k+1)|)  & \mbox{if} & k < d-2, \smallskip \\
\bo( \frac{1}{\sqrt{N_d}}|\sing_1| + |\sing_2|)  & \mbox{if} & k = d-2. 
\end{array}
\right. 
\] 
This establishes (\ref{sing4_est}):
\[ |\sing_4| \leq  \frac{c_2}{\sqrt{d}}  (|\sing_1| + |\sing_2|) \] 
for some constant $c_2 > 0$.
\end{enumerate}

Thus we have 
\begin{eqnarray*}
\pd & = & \frac{|\sing|}{|\reg| + |\sing|}  \\
& = & \frac{(|\sing_1| + |\sing_2|)(1 + \so(1))}{|\reg| + (|\sing_1| + |\sing_2|)(1 + \so(1))}  \\
& = & \frac{ \left( \frac{d-1}{2} + \frac{\ed}{d} \right)  }
{ \left( \frac{d-1}{2} + \frac{\ed}{d}  
 +  \frac{2^d - \ed}{d} \right)  }(1 + \so(1)) \\
& = &  \left( \frac{1}{2^d} \ed + \frac{d^2}{2^{d+1}} \right) (1 + \so(1)) 
\end{eqnarray*}

This concludes the proof of Theorem \ref{th1}. 
\begin{flushright} $\Box$ \end{flushright}

\section{Experiments in small dimensions}
\label{sec:experiments}

Complete enumeration of the $0/1$-matrices of size $d\times d$
is feasible up to dimension 7
(see \cite{Zie99}), while hyperplanes were enumerated up to dimension
8 (see Aichholzer \& Aurenhammer \cite{aichholzer96:_class}).  For
some higher dimensions we generated 25,000,000 random matrices and
determined an experimental probability $P_x(d)$ that a random matrix
is singular.  The significance of these numbers is limited
for high dimensions (we found very few singular matrices and 25
million is tiny compared to the number of $0/1$-matrices),
but since the number of singular matrices is sharply concentrated
around the expected value the results should still
be close to the real values. Up to dimension~17 $P_x(d)$ decreases at a slower rate
than the natural lower bound $d^2 2^{-d}$ while in higher dimensions
$P_x(d)$ seems to approach this bound.
\begin{small}
\[
\begin{array}{r|rrrrrrr}
d & \mbox{matrices} & \mbox{singular} & P_x(d) & 
\frac{d^2}{2^d} & P_x(d) 2^d d^{-2} \\
\hline 
1  & 2^1    &        1 & 0.5000000 & 0.500000 & 1.000 \\
2  & 2^4    &       10 & 0.6250000 &1.000000 & 0.625 \\
3  & 2^9    &      338 & 0.6601562 &1.125000 & 0.587 \\
4  & 2^{16} &    42976 & 0.6557617 &1.000000 & 0.666 \\
5  & 2^{25} & 21040112 & 0.6270442 & 0.781250 & 0.803 \\
6  & 2^{36} & \approx 3.98 \cdot 10^{10} 
                         & 0.5803721 & 0.562500 & 1.032 \\
7  & 2^{49}    & \approx 2.92 \cdot 10^{14} 
                         & 0.5197696 & 0.382812 & 1.358 \\
8  & 25000000 & 11230864 & 0.4492346 & 0.250000 & 1.797 \\
9  & 25000000 &  9331895 & 0.3732758 & 0.158203 & 2.359 \\
10 & 25000000 &  7430305 & 0.2972122 & 0.0976562 & 3.043 \\
11 & 25000000 &  5657196 & 0.2262879 & 0.0590820 & 3.830 \\
12 & 25000000 &  4108304 & 0.1643321 & 0.0351562 & 4.674 \\
13 & 25000000 &  2837245 & 0.1134898 & 0.0206299 & 5.501 \\
14 & 25000000 &  1868850 &  0.0747540 & 0.0119629 & 6.249 \\
15 & 25000000 &  1175425 &  0.0470170 & 0.0068665 & 6.847 \\
16 & 25000000 &   707571 &  0.0283028 & 0.0039062 & 7.246 \\
17 & 25000000 &   407077 &  0.0162831 & 0.0022049 & 7.385 \\
18 & 25000000 &   225820 &  0.0090328 & 0.0012360 & 7.308 \\
19 & 25000000 &   121157 &  0.0048463 & 0.0006886 & 7.038 \\
20 & 25000000 &    62500 &  0.0025000 & 0.0003815 & 6.554 \\
21 & 25000000 &    31779 &  0.0012712 & 0.0002103 & 6.045 \\
22 & 25000000 &    15393 &  0.0006157 & 0.0001154 & 5.336 \\
23 & 25000000 &     7383 &  0.0002953 & 0.0000631 & 4.683 \\
24 & 25000000 &     3515 &  0.0001406 & 0.0000343 & 4.095 \\
25 & 25000000 &     1722 &  0.0000689 & 0.0000186 & 3.698 \\
26 & 25000000 &      736 &  0.0000294 & 0.0000101 & 2.923 \\
27 & 25000000 &      345 &  0.0000138 & 0.0000054 & 2.541 \\
28 & 25000000 &      164 &  0.0000066 & 0.0000029 & 2.246 \\
29 & 25000000 &       81 &  0.0000032 & 0.0000016 & 2.068 \\
30 & 25000000 &       37 &  0.0000015 & 0.0000008 & 1.766 \\
\end{array}
\]
\end{small}

\vspace{20mm}
\textbf{Note added in proof:}\\
  Recently, T. Tao and V. H. Vu \cite{TV04a} 
  have significantly improved the upper bound on $P_s(d)$,
  by proving that $P_s(d)=(\frac34+o(d))^d$. 

  Furthermore, M. \v{Z}ivkovi\'c has recently 
  computed the number of singular $0/1$-matrices of
  size $8\times8$ exactly \cite{zivkovic-2005}.
  From this we get that $P_s(8)=0.4492003726$,
  so our estimate $P_x(8)=0.4492346$ wasn't bad.

\bibliography{BasicBib}
\end{document}